\documentclass{amsart}
\usepackage[english]{babel}
\usepackage{mathrsfs}
\usepackage{epsfig,afterpage}
\usepackage{psfrag}
\usepackage{graphicx}
\usepackage{amsmath}
\usepackage{amssymb}
\usepackage{amsfonts}
\usepackage{multirow}

\newcommand{\const}{{\rm const}}
\newcommand{\bR}{{\mathbb R}}
\newcommand{\bN}{{\mathbb N}}
\newcommand{\bZ}{{\mathbb Z}}

\let\<\langle
\let\>\rangle

\newtheorem{theorem}{Theorem}
\newtheorem{lemma}{Lemma}
\newtheorem{remark}{Remark}
\newtheorem{corollary}{Corollary}
\newtheorem{example}{Example}
\newtheorem{definition}{Definition}

\numberwithin{equation}{section}  

\begin{document}

\title[Singular differential equations]{On solutions of singular differential equations of the second order}

\author[A.O. Remizov]{A.O. Remizov}
\address{
Moscow Institute of Physics and Technology,
Inststitutskii per. 9, 141700 Dolgoprudny, Russia}
\email{alexey-remizov@yandex.ru}




\begin{abstract}
We study the behaviour of solutions of second-order differential equations at their singular points,
where the coefficient of the second-order derivative vanishes. In particular, we consider solutions issuing from
a singular point without definite tangential direction. Great attention is paid to second-order differential equations, whose
right-hand sides are cubic polynomials by the first-order derivative.
\end{abstract}

\maketitle

\begin{flushright}
\it {Dedicated to the memory of A.\,F.~Filippov}
\end{flushright}

\section{Introduction}

A.\,F.~Filippov showed in \cite{Fil} that a system of ordinary differential equations
\begin{equation}
F(x,y,p) = 0, \ \ p=dy/dx,
\label{1}
\end{equation}
where $x \in \bR^1$, $y \in \bR^n$ and $F: \bR^{2n+1} \to \bR^n$ is a smooth\footnote{
By smooth we mean $C^\infty$, if the otherwise is not stated.
}
vector function, may have solutions without definite tangential direction at a certain point.
This motivates the following terminology.

A solution $y(x)$ of system \eqref{1} is called {\it oscillating} at a point $(x_0,y_0)$, if $y(x)$ is a vector function
differentiable on the interval $(x_0, x_0+\delta)$ or $(x_0-\delta,x_0)$, $\delta>0$, such that
$y(x) \to y_0$ but $y'(x)$ has neither finite not infinite limit as $x \to x_0$.
A solution $y(x)$ is called {\it proper} at a point $(x_0,y_0)$, if $y'(x)$ has limit (finite or infinite) as $x \to x_0$.

Let $J^1$ be the 1-jet space of smooth functions $y(x)$, with the coordinates $(x,y,p)$.
Filippov showed that if for a point $T_0 \in J^1$ there exists $T_0' \in J^1$ such that $x_0=x_0'$, $y_0=y_0'$, $p_0 \neq p_0'$,
and the matrix $F_p$ degenerates at $T_0'$, then besides a unique proper solution passing through $(x_0,y_0)$ with the tangential direction $p_0$,
system \eqref{1} may have oscillating solutions, which pass through $(x_0,y_0)$ and even have definite tangential direction $p_0$ at $(x_0,y_0)$.

The following example is taken from~\cite{Fil}. Consider solutions of the system
\begin{equation}
\left\{
\begin{aligned}
& p_1(1-p_1^2-p_2^2) + 8xy_1 +4y_2 = 0, \\
& p_2(1-p_1^2-p_2^2) + 8xy_2 -4y_1 = 0,
\end{aligned}
\right.
\label{2}
\end{equation}
that pass through the origin with the tangential direction $p_i=0$, that is, satisfy the conditions $y_i(0)=y'_i(0)=0$, $i=1,2$.
The matrix $F_p$ is non-degenerate at the point $T_0=0$, and system \eqref{2} has a unique proper solution $y_i(x) \equiv 0$, which
obviously satisfies the required conditions. However, system \eqref{2} has an infinite number of oscillating solutions given by the formula
$$
y_1(x) = x^2 \cos (x^{-1}+c), \ \
y_2(x) = x^2 \sin (x^{-1}+c), \ \
c=\const,
$$
for $x \neq 0$, and $y_i(0)=0$.
The both derivatives $y_i'(x)$ are zero at $x=0$, but the limits $y_i'(x)$ as $x \to 0$ do not exist.

On the other hand, Filippov showed that systems with oscillating solutions are a sort of exception.
Given a point $q_{0}=(x_{0},y_{0}) \in \bR^{n+1}$, we define the set
$$
Q(q_{0}) = \{ T=(q_0,p) \in J^1 : F(T)=0, \ \det (F_p(T))=0 \}.
$$

\begin{theorem}[\cite{Fil}]
\label{T0}
Assume that at least one of the following conditions holds true:

1. $p_0 \notin {\rm co}\,Q(q_{0})$, where ${\rm co}$ denotes the convex hull,

2. the set $Q(q_{0})$ is at most countable (i.e., countable, or finite, or empty).

Then system \eqref{1} has no solution oscillating at $q_{0}$ with tangential direction $p_0$ at $q_{0}$.
Generically, the condition~2 holds true at all points, and system \eqref{1} has no oscillating solutions.
\end{theorem}

Although Theorem~\ref{T0} has a very general character, it gives only a trivial result for systems
\begin{equation}
A(x,y)\,p=b(x,y), \ \ p=dy/dx,
\label{3}
\end{equation}
where $A$ is a $n \times n$ matrix depending on $(x,y)$, $b$ is a vector function.
Theorem~\ref{T0} guarantees the absence of oscillating solution at those points where $\det A \neq 0$,
which obviously follows from the fact that then system \eqref{3} is locally equivalent to $p = A^{-1}(x,y)b(x,y)$.

System \eqref{3} are of great interest themselves (see, e.g., \cite{ZhS} and the references therein), but in this paper
we restrict ourselves with a partial case of \eqref{3}. Namely, we shall consider the second-order equations
\begin{equation}
\Delta(x,y) \frac{dp}{dx} = M(x,y,p),  \ \ p=dy/dx,
\label{4}
\end{equation}
which are equivalent to systems \eqref{3} of a special type.

A great attention is paid to the case that the right-hand side of equation \eqref{4}
is a cubic polynomial in $p$, that is,
\begin{equation}
M(x,y,p) = \sum_{i=0}^3 \mu_i(x,y)p^i.
\label{5}
\end{equation}
This attention to equations \eqref{4}, \eqref{5} is motivated by their role in the description of various geometric structures,
for instance, the affine connection and the projective connection.
Equations of this class were studies by Sophus Lie, A.~Tresse, J.~Liouville, E.~Cartan, to name a few.
(See, e.g., \cite{AA, Cartan, Tresse, Yu} and the references therein.)
For instance, geodesics in the Levi--Civita connection generated by the metric tensor
\begin{equation}
ds^2 = a(x,y)\,dx^2 + 2b(x,y)\,dxdy + c(x,y)\,dy^2
\label{6}
\end{equation}
are solutions of the equation in the form \eqref{4}, \eqref{5}, where
\begin{equation}
\begin{aligned}
\Delta &= ac-b^2,  \\
\mu_0 &= a(a_y-2b_x) +a_xb, \\
\mu_1 &= b(3a_y-2b_x) + a_xc - 2ac_x, \\
\mu_2 &= b(2b_y-3c_x) + 2a_yc - ac_y, \\
\mu_3 &= c(2b_y-c_x) - bc_y.
\end{aligned}
\label{7}
\end{equation}

Here the coefficient $\Delta$ vanishes at point that the quadratic form \eqref{6} degenerates.
In Riemannian geometry, such is not the case, since the inequality $\Delta>0$ holds true everywhere.
However, degenerate quadratic forms generically appear on surfaces embedded into pseudo-Euclidean spaces:
the metric induced on a surface into pseudo-Euclidean space degenerates at those points where the surface tangents
the light cone of the ambient space. Singularities of the geodesic equation appearing at degenerate points of the metric \eqref{6} are
studied in a recent series of papers \cite{PR2019, R2009, RT2016}. See also the survey \cite{PR2018}.

It should be remarked that in all papers mentioned above (as well as in other works known to us)
only proper geodesics were considered, while the possibility of oscillating geodesics is not studied yet.
In the present paper, we shall fill this gap.

The paper is organized as follows.
In the next section, we establish a sufficient condition for equation \eqref{4}, \eqref{5} to have no oscillating solutions at a given point.
This implies that a generic equation of this type has no oscillating solutions at all.
The closing section of the paper is devoted to the behaviour of proper solutions of equations \eqref{4}, \eqref{5}
at its generic singular points.

\section{Oscillating solutions}

Consider the differential equation
\begin{equation}
\Delta(x,y) \frac{dp}{dx} = M(x,y,p),  \ \ p=dy/dx,
\label{8}
\end{equation}
where $\Delta(x,y)$, $M(x,y,p)$ are smooth functions.

A point $q_0=(x_0,y_0)$ is called a {\it singular point} of equation \eqref{8}, if $\Delta(x_0,y_0)=0$.
Denote the locus of singular points by $\Gamma$.
Generically, $\Gamma$ is a curve on the $(x,y)$-plane, but we do not assume it in what follows.
When dealing with singular points, it makes sense to refine the notion of solutions.

Let $q_{0} = (x_{0},y_{0}) \in \Gamma$.
A solution of equation \eqref{8} {\it issuing} from the point $q_{0}$ is a function $y(x)$ such that

1) \ $y(x)$ is continuous on the segment $I_{\varepsilon}$ with endpoints $x_{0}$, $x_{0}+\varepsilon$, $\varepsilon>0$ or $\varepsilon<0$,

2) \ $y(x)$ is differentiable and it satisfies \eqref{8} at all inner points of $I_{\varepsilon}$,

3) \ $y(x_{0})=y_{0}$ and the graph of $y(x)$ has no common points with $\Gamma$ except for $q_0$.

\begin{figure}[!htp]
\centering
\includegraphics[scale=0.85]{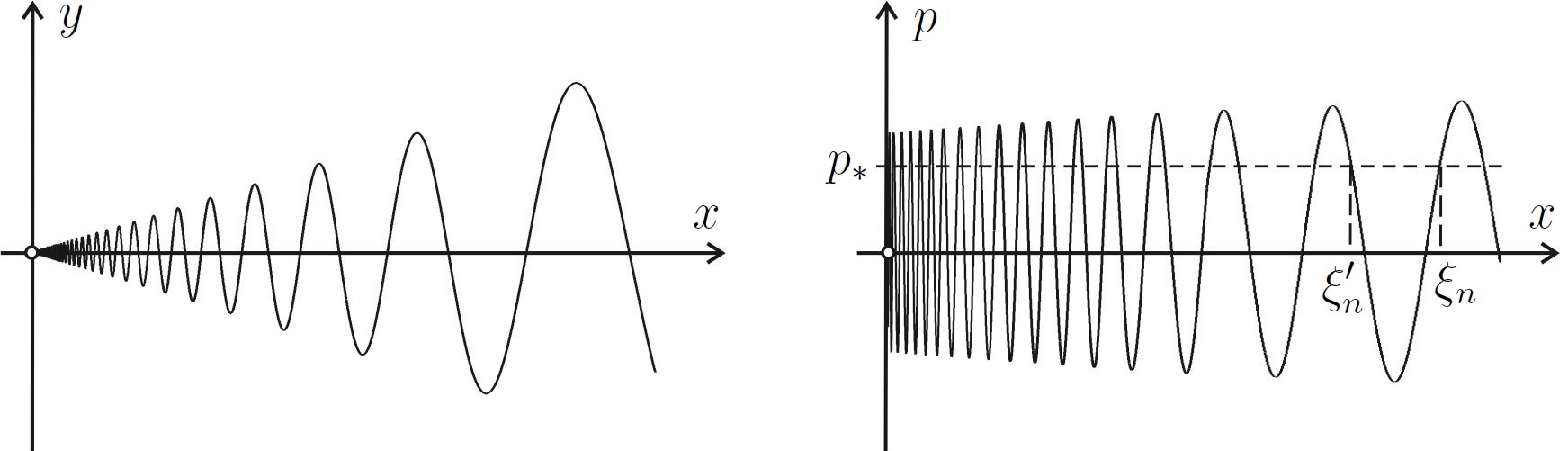}
\caption{An oscillating solution issuing from the origin (on the left) and its derivative (on the right).
}
\label{pict1}
\end{figure}

If in addition to the above conditions, the derivative $y'(x)$ has a (finite or infinite) limit as $x \to x_{0}$,
the solution $y(x)$ is called {\it proper}. Otherwise the solution $y(x)$ is called {\it oscillating} at~$q_{0}$.
See Fig.~\ref{pict1}.

Finally, we shall give one more definition.
A solution of equation \eqref{8} {\it passing through} the point $q_{0}$ is a function $y(x)$ differentiable and satisfying \eqref{8} at all points
of an open interval $I$ that contains $x_0$ such that $y(x_{0})=y_{0}$ the graph of $y(x)$, $x \in I$, does not intersects $\Gamma$ except for the point $q_0$. Obviously, any solutions passing through $q_{0}$ is the union of two solutions issuing from $q_{0}$, but the inverse is not true.

The advantage of the given definitions becomes clear from the following simple example.

\begin{example}
\label{Ex1}
{\rm
The parabolas $y = \alpha x^2$, $\alpha =\const$, are integral curves of the second-order equation $2y {dp}/{dx} = p^2$
passing through the origin.
However, the graphs of functions defined by the formula $y = \alpha_1 x^2$ for $x < 0$ and
$y = \alpha_2 x^2$ for $x \ge 0$ are also integral curves, and all such functions are solution of the given equation
passing through the origin.
To avoid such ambiguity, it is sufficient to use the notion of a solution issuing from a point defined above.
Then solutions issuing from the origin are the branches of the parabolas $y = \alpha x^2$, $\alpha =\const$, and only them.
}
\end{example}

In the previous example, all solutions are proper.
Now we give two examples of equations \eqref{8} that possess oscillating solutions.

\begin{example}
\label{Ex2}
{\rm
The equation $x^4 {dp}/{dx} = 2x^3p - (2x^2+1)y$ has the family of solutions defined by the formula
\begin{equation}
y = x^2  (\alpha \cos x^{-1} + \beta \sin x^{-1}) \ \ \textrm{for} \ \ x \neq 0,  \ \ \alpha, \beta = \const,
\label{X9}
\end{equation}
and zero for $x=0$.
These functions are differentiable on the whole real axis, but their derivatives are not continuous at the origin except for $\alpha =\beta = 0$.
According to the given definitions, formula \eqref{X9} gives a family of solutions passing through the origin with the tangential direction $p=0$.
For $\alpha =\beta = 0$, we have the regular solution $y=0$, while for all others the solutions \eqref{X9} are oscillating at the origin.
}
\end{example}

\begin{example}
\label{Ex3}
{\rm
The equation $x^2 {dp}/{dx} = xp-2y$  has the family of solutions
\begin{equation}
y= x  (\alpha \cos \ln |x| + \beta \sin \ln |x|), \ \ \alpha, \beta = \const,
\label{X10}
\end{equation}
issuing from the origin.
Except for $\alpha =\beta = 0$, all these solutions are oscillating at the origin and (unlike the previous example)
they  have no definite tangential directions at the origin.

Moreover, using formula \eqref{X10}, one can construct an equation of the form \eqref{8} with oscillating solution,
whose right-hand side is a polynomial of arbitrary degree in $p$.
For instance, one can see that the function \eqref{X10} with $\alpha=\beta=1$
is a solution of the first-order equations $f_i(x,y,p)=0$, where
\begin{equation*}
\begin{aligned}
f_2 &= (xp)^2 - 2xyp + 2(y^2-x^2), \\
f_3 &= (xp)^3 + y(xp)^2 - 2xp (x^2+2y^2) + 6y(y^2-x^2).
\end{aligned}
\end{equation*}
Therefore, $y= x(\cos \ln |x| + \sin \ln |x|)$ is an oscillating solution of the second-order equations
$$
x^2 \frac{dp}{dx} = xp-2y + f_i(x,y,p), \ \ p=dy/dx,
$$
where $f_i$ is a polynomial in $p$ of the degree~$i=2,3$.
}
\end{example}

We presented several examples demonstrating the existence of oscillating solutions.
The following theorem shows that oscillating solutions do not appear generically.

\begin{theorem}
\label{T1}
Let $q_{0} \in \Gamma$ and $M(q_{0},p)$ is an analytic function not identically zero.
Then equation \eqref{8} has no oscillating solutions issuing from~$q_{0}$.
\end{theorem}

{\bf Proof.}
Assume that the above conditions hold true, but equation \eqref{8} has an oscillating solution $y(x)$ that issues from $q_0$
without definite tangential direction. Without loss of generality, assume that $y(x)$ is defined on the segment
$I_{\varepsilon} = [x_0, x_0+\varepsilon]$ and $\Delta(x,y(x))>0$ at all its inner points.

From the absence of the limit $y'(x)$ as $x \to x_{0}$ it follows that there exist two sequences
$x_n' \to x_{0}+0$ and $x_n'' \to x_{0}+0$ such that $p(x_n') \to p'$  and $p(x_n'') \to p''$, $p' \neq p''$.
For definiteness, assume that $p' < p''$. Since the function $y'(x)$ is continuous on the interval $x_{0} <x< x_{0}+\varepsilon$,
for every $p_* \in (p', p'')$ there exist two sequences $\xi_n \to x_{0}+0$ and $\xi_n' \to x_{0}+0$ such that
\begin{equation}
\lim_{n \to \infty} p(\xi_n) = \lim_{n \to \infty} p(\xi_n') = p_*, \ \
\frac{dp}{dx}(\xi_n)>0, \ \ \frac{dp}{dx}(\xi_n')<0.
\label{X11}
\end{equation}
See Fig.~\ref{pict1} (right).
Substituting the solution $y(x)$  into \eqref{8}, at the points $x=\xi_n$ and $x=\xi_n'$ we have the equalities
\begin{equation}
\Delta(\xi_n,y(\xi_n)) \frac{dp}{dx}(\xi_n) = M(\xi_n,y(\xi_n),p(\xi_n)),
\label{X12-1}
\end{equation}
\begin{equation}
\Delta(\xi_n',y(\xi_n')) \frac{dp}{dx}(\xi_n') = M(\xi_n',y(\xi_n'),p(\xi_n')).
\label{X12-2}
\end{equation}

The right-hand sides of both equalities \eqref{X12-1}, \eqref{X12-2} have the finite limit $M(x_{0},y_{0},p_*)$, whence their
left-hand sides have the same limit $M(x_{0},y_{0},p_*)$. On the other hand,
$\Delta(\xi_n,y(\xi_n))$ and $\Delta(\xi_n',y(\xi_n'))$ are positive for every $n$, and
from \eqref{X11} it follows that the left-hand sides of equalities \eqref{X12-1}, \eqref{X12-2} have different signs.
Therefore, their limit $M(x_{0},y_{0},p_*)=0$.

Thus, we proved that $M(q_{0},p_*)=0$ for every $p_* \in (p', p'')$.
Since the function $M(q_{0},p)$ is analytic in $p$, this implies $M(q_{0},p)=0$ for every~$p$.
This completes the proof.

\medskip

\begin{corollary}
\label{T2}
Let $q_0$ be a singular point of equation \eqref{4}, \eqref{5}.
If the coefficients $\mu_0, \ldots, \mu_3$ of the cubic monomial $M$ do not vanish at $q_{0}$ simultaneously,
there are no oscillating solutions issuing from $q_{0}$.
Moreover, in the case of the geodesic equation,
$d\Delta(q_{0}) = 0$ is a necessary condition for the existence of an oscillating solution issuing from $q_{0}$.
\end{corollary}

{\bf Proof.}
The first statement is obvious, let us prove the second one.

The geodesic equation in metric \eqref{6} has the form \eqref{4}, \eqref{5}, where $\Delta$ and $\mu_0, \ldots, \mu_3$ are expressed through $a,b,c$
by formula~\eqref{7}. By Theorem~\ref{T1}, for the non-existence of solutions oscillating at $q_{0} \in \Gamma$, it is sufficient to prove
that if all coefficients $\mu_0, \ldots, \mu_3$ simultaneously vanish at $q_{0}$, then $d\Delta (q_{0}) = 0$.

Let us assume on the contrary that $\mu_i(q_{0})=0$ for all $i$ and $d\Delta(q_{0})\neq 0$.
To simplify the calculations, we choose local coordinates centered at $q_{0} \in \Gamma$ such that $b(q_{0})=0$,
which is always possible to attain by an appropriate linear transformation.
Then from the equalities $\Delta(q_{0})= 0$, $d\Delta(q_{0})\neq 0$ it follows that $|a(q_{0})| + |c(q_{0})| \neq 0$.
Without loss of generality one can assume that $a(q_{0}) \neq 0$ and $c(q_{0}) = 0$.
Then from the equalities $\mu_i(q_{0})=0$ and \eqref{7} it follows $c_x(q_{0}) = c_y(q_{0})=0$.
Finally, from the equalities
$$
b(q_{0}) = c(q_{0}) = c_x(q_{0}) = c_y(q_{0})=0
$$
it follows $d\Delta (q_{0})=0$.
The obtained contradiction completes the proof.

\section{Singularities of second-order equations cubic in the first-order derivative}

In this section, we shall examine generic singularities of equation \eqref{4}, whose right-hand side has the form \eqref{5}.
Here we assume that $\Delta(x,y)$ and $\mu_i(x,y)$ are smooth functions not connected with each other, in contrast with the geodesic equation, whose
coefficients are expressed through three functions $a,b,c$ (formula~\eqref{7}), and therefore, they are connected with some
functional relations. According to the genericity assumption, we shall assume that the coefficients $\mu_0, \ldots, \mu_3$ of the cubic monomial \eqref{5}
do not vanish simultaneously and the set of singular points $\Gamma$ is a regular curve, which locally separates the $(x,y)$-plane
into two domains.

Under the above conditions, equation \eqref{4}, \eqref{5} has no oscillating solutions (Theorem~\ref{T1}).
Our purpose is to study its proper solutions issuing from generic singular points.
Further, we shall omit the adjective {\it proper}.
We start with a brief survey of the basic facts established in the papers \cite{PR2018} -- \cite{RT2016}
for geodesic equations that are valid for all equations of the form \eqref{4}, \eqref{5}, not necessarily geodesic.

Denote by $\pi$ the standard projection acting from the space $J^1$ to the $(x,y)$-plane along to the $p$-direction, which we call {\it vertical},
that is, $\pi: (x,y,p) \mapsto (x,y)$. Equation \eqref{4}, \eqref{5} is obviously connected with the vector field
\begin{equation}
\dot x = \Delta(x,y), \ \ \dot y = p\Delta(x,y), \ \ \dot p = M(x,y,p).
\label{X13}
\end{equation}
The $\pi$-projection of any integral curve of the field \eqref{X13} different from a straight vertical line, is a solution of equation \eqref{4}.
Conversely, the Legendrian lift of every solution of \eqref{4} is an integral curve of the field \eqref{X13}.

\begin{lemma}
\label{L0}
For any point $q_{0} \in \Gamma$, the vertical line $\pi^{-1}(q_{0})$ is an integral curve of the field \eqref{X13}.
Moreover, if $p_{0}$ satisfies the condition $M(q_{0},p_{0}) \neq 0$, then $\pi^{-1}(q_{0})$
is the only integral curve of the field \eqref{X13} that passes through the point $(q_{0},p_{0})$.
\end{lemma}

The proof is trivial.

From Lemma \ref{L0} it follows that solutions of equation \eqref{4}, \eqref{5} can issue from a singular point $q_{0} \in \Gamma$
only in so-called {\it admissible} tangential directions $p$, which correspond to real roots of the cubic polynomial $M(q_{0},p)$.
Here we shall understand the term {\it cubic polynomial} in the broad sense: the higher coefficients $\mu_3$ can vanish at some points
or even identically. Under such a convention, the class of equations \eqref{4}, \eqref{5} is invariant with respect to changes
of the variables $x,y$. It includes, in particular, the Bessel equation and the Gaussian hypergeometric equation.

The tangential direction $p=\infty$ is admissible at $q_{0}$, if and only if $\mu_3(q_{0})=0$.
Indeed, interchanging the variables $x,y$, we obtain the equation
$$
\Delta(y,x) \frac{dp}{dx} = -M^*(y,x,p) = - \sum_{i=0}^{3} \mu_i(x,y) p^{3-i},
$$
where $M^*$ is the reciprocal polynomial, and the direction $p=\infty$ becomes $p=0$.
The equality $\mu_3(q_{0})=0$ is equivalent to the condition that $p=0$ is a root of $M^*(q_{0},p)$.
Therefore, we shall consider $p=\infty$ as a root of the cubic polynomial $M(q_{0},p)$ of multiplicity $k$ iff
$p=0$ as a root of the cubic polynomial $M^*(q_{0},p)$ of multiplicity~$k$.

\begin{definition}
We shall call a singular point of equation \eqref{4}, \eqref{5} {\it generic}, if
$d\Delta(q_0) \neq 0$,
all roots of the cubic polynomial $M(q_{0},p)$ are prime and the corresponding admissible directions are
transversal to the curve $\Gamma$.
\end{definition}

In a neighborhood of every generic singular point $q_{0}$, there exist local coordinates such that $\mu_3(q_{0}) \neq 0$,
and consequently, all admissible directions are finite. In what follows, we shall assume that this condition is satisfied.
Then the polynomial $M(q_{0},p)$ has either a unique real $p_0$ or three different real roots $p_0, p_1, p_2$.

By Lemma \ref{L0}, all solutions of equation \eqref{4}, \eqref{5} that issue from a generic point $q_{0} \in \Gamma$
are $\pi$-projections of integral curves of the field \eqref{X13} that tend to its singular point $(q_{0},p)$ as
time tends to infinity (positive or negative). Singular points of the field \eqref{X13} are given by the equations
\begin{equation}
\Delta(x,y) = 0, \ \  M(x,y,p)=0,
\label{X14}
\end{equation}
which define a single curve or three non-intersecting curves, each of them corresponds to one of the roots $p_0, p_1, p_2$.
The spectrum $\Sigma$ of the linear part of \eqref{X13} at its singular point $(q_{0},p)$ has the form
\begin{equation}
\Sigma = (0, \lambda_1, \lambda_2), \ \ \lambda_1=\Delta_x + p\Delta_y, \ \ \lambda_2=M_p.
\label{X15}
\end{equation}

From the condition of genericity (see the above definition) it follows that $\lambda_{1,2} \neq 0$.
It is worth observing that $\partial_p$ is the eigenvector with the eigenvalue $\lambda_2=M_p$, while
the eigenvector with the eigenvalue $\lambda_1$ is not vertical.
Denote
\begin{equation}
\lambda = \lambda_2/\lambda_1.
\label{X16}
\end{equation}

Recall several basic definitions concerning local equivalence of vector fields (see, e.g., \cite{AI}).

Let $V_1$ and $V_2$ be the germs of two smooth vector fields in $\bR^n$.
The germs $V_1$ and $V_2$ are called $C^k$-smoothly (resp. topologically) equivalent, if there exists a $C^k$-diffeomorphism (resp.  homeomorphism)
$h \colon \bR^n \to \bR^n$ that conjugates their phase flows $g^t_1$ and $g^t_2$, i.e., $h \circ g^t_1 \equiv g^t_2 \circ h$.
Here $k \ge 1$ is an integer number (finite-smooth equivalence) or $\infty$ (infinite-smooth equivalence).
The germs $V_1$ and $V_2$ are called {\it orbitally} $C^k$-smoothly (resp. topologically) equivalent,
if there exists a $C^k$-diffeomorphism (resp.  homeomorphism) $h \colon \bR^n \to \bR^n$
that conjugates their integral curves (orbits of their phase flows).
In other words, $V_1$ and $V_2$ are {orbitally} $C^k$-smoothly (resp. topologically) equivalent,
if there exists a non-vanishing scalar function $\varphi$ such that $V_1$ and $\varphi V_2$ are $C^k$-smoothly (resp. topologically) equivalent.

The second definition (orbital equivalence) slightly differs from the generally accepted definition of the orbital equivalence,
where coincidence of the orientation of integral curves is also required.
In fact, our definition is naturally related to directions fields, whose integral curves
do not have an orientation a priori.

Consider the germ of a smooth vector field $V$ in $\bR^3$, whose component belong to the ideal $I$
(in the ring of smooth functions) generated by two of them.
In appropriate local coordinates, such a field has the form
\begin{equation}
\dot \xi = v, \ \
\dot \eta = w, \ \
\dot \zeta = \alpha v + \beta w,
\label{X17}
\end{equation}
where $v,w$ and $\alpha, \beta$ are arbitrary smooth functions on the variables $\xi, \eta, \zeta$.
The set of singular points of \eqref{X17} is given by two equations $v=w=0$ and the spectrum of the linear part of \eqref{X17} at every singular point
has the form $\Sigma = (\lambda_1, \lambda_2, 0)$. For instance, vector field \eqref{X13} has the form \eqref{X17}, in this case
the ideal $I$ is generated by the components $\Delta$ and $M$.

Let $T_0$ be a singular point of \eqref{X17} such that the eigenvalues $\lambda_{1,2}$ are real and non-zero.
Then in a neighborhood of $T_0$, the set of singular points is the center manifold of \eqref{X17}, and the reduction principle
(see \cite{AI, HPS}) leads to the following result:

\begin{lemma}
\label{L1}
The germ of vector field \eqref{X17} at $T_0$ is orbitally topologically equi\-va\-lent to
\begin{equation}
\dot \xi = \xi, \ \ \dot \eta = \pm \eta, \ \ \dot \zeta = 0.
\label{X18}
\end{equation}
For vector field \eqref{X13}, the sign $\pm$ coincides with the sign of $\lambda$ defined in \eqref{X16}.
\end{lemma}

\medskip

Lemma~\ref{L1} shows that vector field \eqref{X17} has a local invariant foliation that is topologically equivalent to $\zeta = \const$.
Under certain additional conditions, this invariant foliation can be brought to the form $\zeta = \const$ by means of a local diffeomorphism,
finitely or infinitely smooth. The main obstacle for the existence of such a diffeomorphism is the resonance relations
\begin{equation}
p \lambda_1 + q\lambda_2 = 0, \ \ p,q \in \bZ_+, \ \  p+q \ge 1,
\label{X19}
\end{equation}
the integer number $p+q$ is called the order of resonance \eqref{X19}.
For instance, using the results of \cite{Sam1}, it is not hard to prove the following statement.

Following \cite{Sam1}, for every integer $k \ge 1$ we define the number $N(k) = 2 \left[ (2k+1) m_1/m_2 \right] + 2$, where
\begin{equation*}
m_1 = \max |\lambda_{1,2}|, \ \ \, m_2 = \min |\lambda_{1,2}|,
\end{equation*}
and the square brackets denote the integer part.
If $\lambda_{1,2}$ have no resonances \eqref{X19} up to the order $N(k)$,
then there exists a local $C^k$-smooth diffeomorphism that brings the germ of \eqref{X17} to the form
\begin{equation}
\dot \xi = X(\xi, \eta, \zeta), \ \ \, \dot \eta = Y(\xi, \eta, \zeta), \ \ \, \dot \zeta = 0,
\label{X21}
\end{equation}
where $X,Y$ are smooth functions that the ideal $I = \<X,Y\> = \<\xi,\eta\>$.

There are examples showing that if $\lambda_{1,2}$ have resonances \eqref{X19}, such diffeomorphism does not exist or at least it has a restricted smoothness. For instance, for a special class of vector fields \eqref{X17}, whose non-zero eigenvalues have resonance \eqref{X19}
with relatively prime $p \ge q$, it is shown that the reduction to \eqref{X21} holds true by mean of $C^{p-1}$-smooth diffeomorphism,
but not $C^{p}$-smooth. See Theorem~5 in \cite{PavRem2020}.

Further we shall deal with the case that $\lambda_{1,2}$ have no resonances \eqref{X19} of all orders.
There are two different cases: $\lambda$ defined in \eqref{X16} is either positive or negative irrational.
In the both cases, for every integer $k \ge 1$ there exists a $C^k$-smooth diffeomorphism
\begin{equation}
f_k: \ U_k \to f_k (U_k),
\label{X22}
\end{equation}
where $U_k$ is a neighborhood of the point $T_0$, that brings \eqref{X17} in $U_k$ to the form \eqref{X21}.
If $\lambda$ is positive, there exists a neighborhood $U$ and a $C^\infty$-smooth diffeomorphism
$f: U \to f(U)$ that brings \eqref{X17} in $U$ to the form \eqref{X21}.
This statement is proved in \cite{Rouss}, a~simpler proof can be obtained by the homotopy method; see \cite{IYa} (Sections 1.2\,--\,1.5).

If $\lambda$ is negative irrational, a $C^\infty$-smooth diffeomorphism $f: U \to f(U)$
that brings \eqref{X17} to the form \eqref{X21} exists only in the exceptionally rare case:
if the value $\lambda$ is the same at all singular points in a neighborhood of~$T_0$.
See, e.g. \cite{Rouss, PavRem2020}. Otherwise, every neighborhood of the point $T_0$ contains a singular point
of the field \eqref{X17} with resonance \eqref{X19}, and there required diffeomorphism $f$ does not exist.

Taking into account the previous reasonings and using results from \cite{IYa, Sam1, Sam2}, one can get the following result.
As before, we put $\lambda = \lambda_2/\lambda_1$. Define the set $\bN^{-1} = \{1/n, \ n \in \bN\}$.

\begin{lemma}
\label{L2}
Assume that $T_0$ is a singular point of vector field \eqref{X17} such that $\lambda$ is either positive or negative irrational.
Then for any integer $k \ge 1$ the germ of vector field \eqref{X17} at the point $T_0$ is $C^k$-smoothly equivalent
to one of the germs at the origin:
\begin{eqnarray}
&&
\dot \xi = a_1(\zeta)\xi, \ \
\dot \eta = a_2(\zeta)\eta, \ \ \phantom{111111111}
\dot \zeta = 0,
\quad \textrm{if} \ \ \lambda \notin \bN \cup \bN^{-1},
\label{X23}\\
&&
\dot \xi = a_1(\zeta)\xi, \ \
\dot \eta = a_2(\zeta)\eta + b_2(\zeta)\xi^n, \ \
\dot \zeta = 0,
\quad \textrm{if} \ \ \lambda_2 = n\lambda_1, \ n>1,
\label{X24} \\
&&
\dot \xi = a_1(\zeta)\xi + b_1(\zeta)\eta^n, \ \
\dot \eta = a_2(\zeta)\eta, \ \
\dot \zeta = 0,
\quad \textrm{if} \ \ \lambda_1 = n\lambda_2, \ n > 1,
\label{X25} \\
&&
\dot \xi = a_1(\zeta)\xi + b_1(\zeta)\eta, \ \
\dot \eta = a_2(\zeta)\eta + b_2(\zeta)\xi, \ \
\dot \zeta = 0,
\quad \textrm{if} \ \ \lambda_1 = \lambda_2,
\label{X255}
\end{eqnarray}
where $n$ is integer, $a_i, b_i$ are smooth functions, $a_i(0)=\lambda_i$.
Moreover, if $\lambda>0$, then the same statement holds true with $k=\infty$.
\end{lemma}

Lemma \ref{L2} states that there exists a $C^k$-smooth diffeomorphism \eqref{X22} bringing the field \eqref{X17} to one of the normal forms
\eqref{X23} -- \eqref{X255}, where $U_k$ is a neighborhood of the point $T_0$ and $f_k(U_k)$ is a neighborhood of the origin,
$f_k(T_0)=0$. All vector fields \eqref{X23} -- \eqref{X255} have the invariant foliation
$$
S_c = \{\zeta = c\}, \ \ c= \const.
$$
This yields the corresponding invariant foliation of the field \eqref{X13} with the leaves $f_k^{-1}(S_c)$.
The restriction of \eqref{X13} to the leaf $f_k^{-1}(S_c)$ is a saddle if $\lambda<0$ or a node if $\lambda>0$.

It is worth observing that \eqref{X23} -- \eqref{X255} are a natural generalization of the Poincar\'e--Dulac normal form
for vector fields \eqref{X17}. The monomial $b_2(\zeta)\xi^n$ in \eqref{X24} and $b_1(\zeta)\eta^n$ in \eqref{X25} are
corollaries of the resonances $\lambda_2 = n\lambda_1$, $\lambda_1 = n\lambda_2$, respectively.
In the case $\lambda_1 = \lambda_2$, the both resonances hold true, whence the normal form \eqref{X255}
contains the both monomials $b_1(\zeta)\eta$, $b_2(\zeta)\xi$.
Generically, none of them can be simplified identically to zero. See, e.g., \cite{Arn71}.

\begin{lemma}
\label{L3}
Consider the germ of vector field \eqref{X13} at its singular point $T_0=(q_0,p_i)$ with non-zero eigenvalues $\lambda_{1,2}$.
If $\lambda_1 = n\lambda_2$ with some integer $n \ge 1$, then the coefficient $b_1(\zeta)$ in normal forms \eqref{X25}, \eqref{X255}
vanishes at $\zeta=0$.
\end{lemma}

\begin{figure}[!htp]
\centering
\includegraphics[width=360pt,height=160pt]{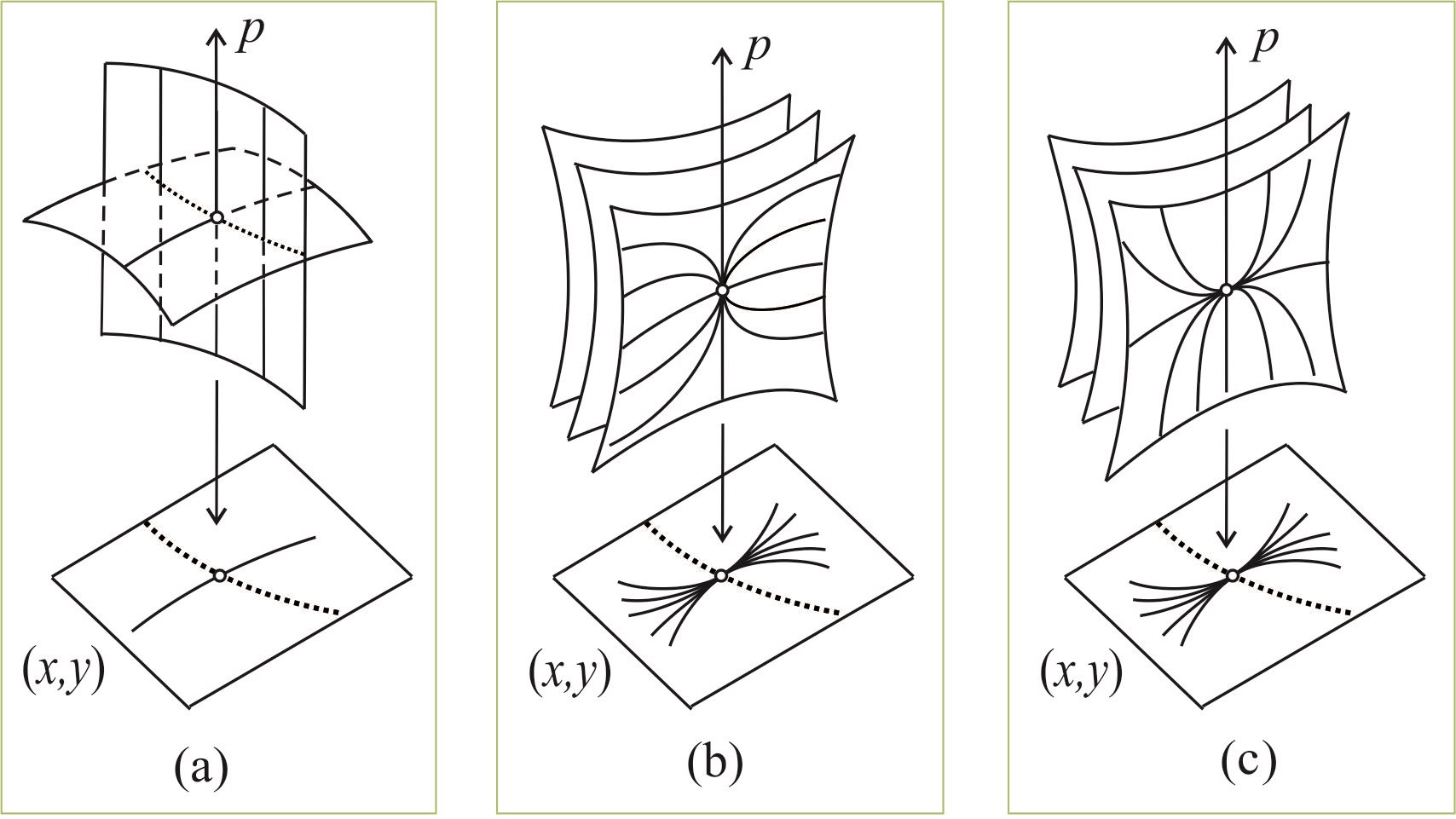}
\caption{
Integral curves of the field \eqref{X13} and their $\pi$-projections passing through a fixed point $q \in \Gamma$.
From left to right: $\lambda<0$ (a), $0<\lambda<1$ (b), $\lambda>1$ (c).
The curve $\Gamma$ is depicted as a dotted line.
}
\label{pict2}
\end{figure}

{\sc Proof}.
Lemma \ref{L2} states that there exists a $C^\infty$-smooth diffeomorphism $f: U \to f(U)$ bringing the field \eqref{X13} to the normal form
\eqref{X25}, where $U$ is a neighborhood of the point $T_0$ and $f(U)$ is a neighborhood of the origin, $f(T_0)=0$.
Consider the restriction of \eqref{X13} to the leaf $f^{-1}(S_0)$ passing through the point $T_0$.
It is a resonant node that has a single integral curve
with the tangential direction corresponding to the eigenvector with $\lambda_1$
and an infinite family of integral curves with the tangential direction $\partial_p$, the eigendirection with $\lambda_2$.
This family contains the vertical straight line $\pi^{-1}(q_0)$. See Fig.~\ref{pict2}~(a).

Consider the case $\lambda_1 = n\lambda_2$ with integer $n > 1$.
Integrating the restriction of the field \eqref{X25} to the invariant plane $S_0$,
one can get the exact formula for the images of integral curves of this family under the diffeomorphism $f$:
\begin{equation}
\xi = \eta^n(c+\varepsilon \ln|\eta|), \ \ \zeta = 0, \ \ c=\const,
\label{X26}
\end{equation}
where $\varepsilon = {b_1(0)}/{a_2(0)}$.
If $\varepsilon=0$, then all curves of the family \eqref{X26} are $C^\infty$-smooth.
If $\varepsilon \neq 0$, they are only $C^{n-1}$-smooth but not $C^{n}$ at the origin.
This contradicts the fact that one of the curves \eqref{X26} is the image of the vertical straight line $\pi^{-1}(q_0)$
under the $C^\infty$-smooth diffeomorphism $f$. Thus, we proved that $\varepsilon=0$, and consequently, $b_1(0)=0$.

The above reasoning is not applicable if $\lambda_1 = \lambda_2$, but in this case the situation is even simpler. The restriction of \eqref{X255} to the invariant plane $S_0$ is a linear vector field, and its matrix can be brought to the Jordan normal form by an appropriate linear transformation.

\medskip

Let $q_{0} \in \Gamma$ be a generic singular point of equation \eqref{4}, \eqref{5}
and $p_i$ be an admissible direction at $q_{0}$, that is, $p_i$ is a real root of the cubic polynomial $M(q_{0},p)$ in $p$.
Then $T_0=(q_0,p_i)$ is the corresponding singular point of vector field \eqref{X13}
with non-zero eigenvalues $\lambda_{1,2}$ defined in \eqref{X15}.
As before, we use the notation $\lambda = \lambda_2/\lambda_1$.

Before formulating this theorem, we remark that from the condition $d\Delta(q_0) \neq 0$ it follows that in a neighborhood of $q_{0}$,
$\Gamma$ is a regular curve that locally splits the $(x,y)$-plane into two open domains.
Denote them by $\Gamma_+$ and $\Gamma_-$.
Then solutions of equation \eqref{4}, \eqref{5} issuing from a generic singular point $q_{0} \in \Gamma$
with tangential direction $p_i$ are described by the following theorem.

\begin{theorem}
\label{T3}
Assume that $\lambda$ is either positive or negative irrational.
Then the following statements hold true:

1. If $\lambda<0$, then there exists a $C^\infty$-smooth solution passing through $q_{0}$
and there are no other solutions issuing from $q_{0}$ with tangential direction $p_i$.

2.
If $\lambda>0$, then there exist an infinite number of solutions issuing from $q_{0}$ with tangential direction $p_i$
in $\Gamma_+$ and an infinite number of solutions issuing from $q_{0}$ with tangential direction $p_i$ in $\Gamma_-$.
There exist local coordinates centered at $q_{0}$ such that the curve $\Gamma$ coincides with the axis $x=0$ and the infinite
family of solutions mentioned above has the form
\begin{eqnarray}
&&
y = F(x, c|x|^{\lambda}),
\quad \textrm{if} \ \ \lambda \notin \bN,
\label{X27}\\
&&
y = F(x, x^n(c+\varepsilon \ln|x|)),
\quad \textrm{if} \ \ \lambda = n \in \bN,
\label{X28}
\end{eqnarray}
where $\varepsilon$ is a real number, $F$ is a smooth function on two variables, $c=\const$.
\end{theorem}

{\bf Proof.}
Recall that solutions of equation \eqref{4}, \eqref{5} issuing from $q_{0}$ with tangential direction $p_i$
are $\pi$-projections of integral curves of the vector field \eqref{X13} that enter its singular point $T_0=(q_0,p_i)$
as time tends to infinity. From Lemmas \ref{L1}\,--\,\ref{L3} it follows that such integral curves of the field \eqref{X13}
and their $\pi$-projections are located as it presented in Fig.~\ref{pict2}
(from the left to right: the cases $\lambda<0$, $0<\lambda<1$, $\lambda>1$).

If $\lambda<0$, then the field \eqref{X13} has two integral curves that enter $T_0$~--
the stable and unstable manifolds of the field, both one-dimensional and $C^\infty$-smooth.
One of them coincides with the vertical straight line $\pi^{-1}(q_{0})$,
and its $\pi$-projection is the point $q_0$.
Another one has a non-vertical tangential direction, and its $\pi$-projection is a unique solution passing
through $q_{0}$ with the tangential direction $p_i$. See Fig.~\ref{pict2}~(a).

If $\lambda>0$, then the field \eqref{X13} has a one-parameter family of integral curves that enter the point $T_0$
as time tends to plus or minus infinity. Taking into account Lemmas \ref{L2}, \ref{L3} and
exploiting the notations used in the proofs of these lemmas,
one can see that this family belongs to the invariant surface $f^{-1}(S_0)$,
where $f: U \to f(U)$ is a diffeomorphism bringing \eqref{X13}
to one of the normal forms \eqref{X23} -- \eqref{X255}.
The restriction of these normal forms to the invariant plane $S_0$ is
\begin{eqnarray}
&&
\dot \xi = \xi, \ \  \dot \eta = \lambda \eta,
\quad \textrm{if} \ \ \lambda \notin \bN,
\label{X29}\\
&&
\dot \xi = \xi, \ \  \dot \eta = \lambda \eta + \varepsilon \xi^{n},
\quad \textrm{if} \ \ \lambda = n \in \bN.
\label{X30}
\end{eqnarray}
Integrating, we get the families of integral curves
\begin{eqnarray}
&&
\{\eta = c |\xi|^{\lambda}, \ \ c=\const\}, \ \ \{\xi = 0\},
\quad \textrm{if} \ \ \lambda \notin \bN,
\label{X31}\\
&&
\{\eta = \xi^n(c+\varepsilon \ln|\xi|), \ \ c=\const\}, \ \ \{\xi = 0\},
\quad \textrm{if} \ \ \lambda = n \in \bN.
\label{X32}
\end{eqnarray}
The $\pi$-projections of the images of the curves \eqref{X31}, \eqref{X32} under the diffeomorphism $f^{-1}$
are infinite families of solutions mentioned in the second statement of the theorem.
See Fig.~\ref{pict2}~(b,\,c). It is important to remark that the family of solutions obtained from the curves \eqref{X31}
contains at least one $C^\infty$-smooth solution, which corresponds either to $\eta=0$ or to $\xi=0$.
On the contrary, if $\varepsilon \neq 0$, then the family of solutions obtained from the curves \eqref{X32}
contains no $C^\infty$-smooth solutions, since the only smooth curve $\xi=0$ corresponds to the vertical straight line $\pi^{-1}(q_0)$.
This fact will be used below.

It remains to prove the existence of local coordinates such that the families of solution obtained from the curves
\eqref{X31}, \eqref{X32}  have the forms \eqref{X27}, \eqref{X28}, respectively.
Let us choose smooth local coordinates on the $(x,y)$-plane centered in $q_{0}$ such that the curve $\Gamma$
coincides with the axis $x=0$ and the admissible direction $p_i=0$.
Moreover, in the case \eqref{X31} one can choose local coordinates such that one of solutions issuing from $q_0=0$
lies on the axis $y=0$.

The invariant surface $f^{-1}(S_0)$ is locally presented as the graph of a smooth function
\begin{equation}
y = F_0(x,p).
\label{X33}
\end{equation}
After multiplication by an appropriate scalar factor, the restriction of the field \eqref{X13} to the invariant surface
$f^{-1}(S_0)$ has the form
\begin{equation}
\dot x = x, \ \ \dot p = p A(x,p) + xB(x),
\label{X34}
\end{equation}
where $A,B$ are smooth functions, $A(0)=\lambda$.
By $f_0$ denote the restriction of the local diffeomorphism $f$ to the surface $f^{-1}(S_0)$, it brings \eqref{X34} to the form \eqref{X29} if
$\lambda \notin \bN$ or to the form \eqref{X30} if $\lambda \in \bN$.
Comparing the vector field \eqref{X34} with \eqref{X29}, \eqref{X30}, one can see that $f_0$ can be chosen in the form
\begin{equation}
x = \xi, \ \ p = \eta \varphi(\xi,\eta) + \xi \psi(\xi,\eta),
\label{X35}
\end{equation}
where $\varphi, \psi$ are smooth functions, $\varphi(0)=1$.
For according to the choice of local coordinates, in the case $\lambda \notin \bN$
the axis $p=0$ is an integral curve of \eqref{X34}, whence $B \equiv 0$ in \eqref{X34} and $\psi \equiv 0$ in~\eqref{X35}.
Thus, if $\lambda \notin \bN$, then from \eqref{X31} and \eqref{X35} with $\psi \equiv 0$
we have the following expression for the derivative of solutions: $p=\eta \varphi(x,\eta)$, where $\eta = c|x|^{\lambda}$.
Substituting it into \eqref{X33}, we have
$$
y = F_0(x,\eta \varphi(x,\eta)) = F(x, \eta) = F(x, c|x|^{\lambda}), \ \ c=\const,
$$
with a smooth function $F$. The remaining integral curve $\{\xi = 0\}$ corresponds to the axis $x=0$ on the $(x,y)$-plane,
and therefore, it does not give a solution.

Similarly, if $\lambda = n \in \bN$, then from formulas \eqref{X32} and \eqref{X35} we obtain the expression
$$
p=\eta \eta \varphi(x,\eta) + x \psi(x,\eta), \ \  \eta = x^n(c+\varepsilon \ln|x|).
$$
Substituting it into \eqref{X33}, we obtain
$$
y = F_0(x,\eta \eta \varphi(x,\eta) + x \psi(x,\eta)) = F(x, \eta) = F(x, x^n(c+\varepsilon \ln|x|)), \ \ c=\const,
$$
with a smooth function $F$.
The proof of theorem is completed.

\begin{remark}
\label{Rem1}
{\rm
Theorem~\ref{T3} is valid
for a single admissible direction $p_i$ that required the only condition that the corresponding value $\lambda$
is either positive or negative irrational, even if other admissible directions at the given singular point do not satisfy this condition.
}
\end{remark}

\begin{example}
\label{Ex4}
{\rm
Integrating the differential equation
\begin{equation}
x \frac{dp}{dx} = \alpha p(p^2-1), \ \ \alpha \neq 0,
\label{X36}
\end{equation}
we get two families
\begin{equation}
p= \pm  \frac{1}{\sqrt{1+c|x|^{2\alpha}}},  \ \ c={\rm const},
\label{X37}
\end{equation}
and the single solution $p=0$, which can be considered as the limit of \eqref{X37} as $c = +\infty$.

The polynomial $M(p)=\alpha p(p^2-1)$ has three roots $p_0=0$, $p_{1,2}=\pm 1$, which correspond to singular points of the vector field
\begin{equation}
\dot x = x, \ \ \dot y = px, \ \ \dot p = \alpha p(p^2-1).
\label{X38}
\end{equation}
At every singular point $x=0$, $p=0$, the value $\lambda = -\alpha$, while
at every singular point $x=0$, $p=\pm 1$, the value $\lambda = 2\alpha$.
This explains the difference between the cases $\alpha > 0$ and $\alpha < 0$.

If $\alpha > 0$, then for every point $q_0 \in \Gamma = \{x=0\}$ there exist
two infinite families of solutions issuing from $q_0$ with tangential directions $p_{1,2}= \pm 1$.
Solutions of the both families are obtained from \eqref{X37} when the constant $c$ runs through the real numbers.
Moreover, there exists the single solution $y=y_0$ passing through the point $q_0$ with the direction $p_0=0$.

If $\alpha < 0$, then there exist two solutions passing through the point $q_0$ with the admissible directions $p_{1,2}$:
$y=y_0 \pm x$, which can be obtained from \eqref{X37} with $c=0$.
Moreover, there exists a one-parameter family of solutions issuing from $q_0$ with the admissible direction $p_0=0$,
which can be obtained from formula \eqref{X37} with all possible $c>0$, including the value $c=+\infty$
(the solution $y=y_0$).

Everything that was said above is in accord with Theorem~\ref{T3}.
All one-parameter families of solutions of equation \eqref{X36} have the form \eqref{X27},
while the form \eqref{X28} with $\varepsilon \neq 0$ does not appear.
The latter property is also in accord with the general theory.
Indeed, the field \eqref{X38} has the form \eqref{X34} with $B \equiv 0$,
which yield $\varepsilon = 0$ in the final form~\eqref{X28}.

In the next example, the logarithm appears.
}
\end{example}

\begin{example}
\label{Ex5}
{\rm
Consider the equation
\begin{equation}
x \frac{dp}{dx} = \alpha p + f(x), \ \ f(x) = \sum_{i\ge 1} f_i x^i, \ \ \alpha \neq 0,
\label{X39}
\end{equation}
where $f(x)$ can be an analytic function or a polynomial of any degree.
From the general viewpoint, here we have the cubic polynomial
$M(p) = \mu_3 p^3 + \mu_2 p^2 + \alpha p + f(x)$, $\mu_2=\mu_3=0$.
Since $f(0)=0$, this yields the admissible direction $p=0$ and the double admissible direction $p=\infty$.
Theorem~\ref{T3} is applicable for studying solutions that issue from singular points of equation \eqref{X39}
with the tangential direction $p=0$ (see Remark~\ref{Rem1}).
It claims that in the case $\alpha<0$ there exists a unique solution of such type, while
in the case $\alpha>0$ there exists an infinite family of such solutions.

Linear equation \eqref{X39} can be explicitly solved, for instance, using the method of variation of parameters.
For simplicity (avoiding modules), we consider it for $x \ge 0$. Then the general solution has the form
$p = A(x) x^{\alpha}$, where $A'(x) x^{\alpha+1} = f(x)$, i.e., $A'(x) = \sum_{}f_i x^{i-\alpha-1}$.
Assume that $\alpha$ in not integer. Then, after integration the latter equality, we get
\begin{equation*}
p = x^{\alpha} A(x) = \sum_{i \ge 1}g_i x^{i} + cx^{\alpha}, \ \ \ g_i = \frac{f_i}{i-\alpha}, \ \ c=\const.
\end{equation*}
Integrating by $x$, we obtain the expression \eqref{X27} for all solutions
that issue from singular points of equation \eqref{X39} with the tangential direction $p=0$.

Now consider the case that $\alpha$ in integer, let us put $\alpha=n$.
Then, after integration the equality $A'(x) = \sum_{}f_i x^{i-\alpha-1}$, we get
\begin{equation*}
p = x^{\alpha} \Bigl( \sum_{i \neq n}g_i x^{i-n} + f_{n} \ln x + c \Bigr) =
\sum_{i \neq n}g_i x^{i} + x^n (c + f_{n} \ln x),
\end{equation*}
where $g_i$, $i \neq n$, are the same as above and $c=\const$.
Integrating by $x$, we obtain the expression \eqref{X28} for all solutions
that issue from singular points with the tangential direction $p=0$.
}
\end{example}

\textbf{Acknowledgements}.
This work is supported
by the Russian Foundation for Basic Research (projects 19-51-50005 and 20-01-00610).
and by the Ministry of Science and Higher Education of the Russian Federation (Goszadaniye N.~075-00337-20-03, project 0714-2020-0005).

\end{document}